\documentclass{amsart}
\usepackage{amssymb}
\usepackage{amsmath}
\newtheorem{theorem}{Theorem}[section]
\newtheorem{lemma}[theorem]{Lemma}
\newtheorem{definition-theorem}[theorem]{Definition-Theorem}
\newtheorem{proposition}[theorem]{Proposition}
\newtheorem{corollary}[theorem]{Corollary}
\newtheorem{definition}[theorem]{Definition}
\newtheorem{example}[theorem]{Example}
\theoremstyle{remark}
\newtheorem{remark}[theorem]{Remark}
\numberwithin{equation}{section}

\newcommand \bth[1] { \begin{theorem}\label{t#1} }
\newcommand \ble[1] { \begin{lemma}\label{l#1} }

\newcommand \bpr[1] { \begin{proposition}\label{p#1} }
\newcommand \bco[1] { \begin{corollary}\label{c#1} }
\newcommand \bde[1] { \begin{definition}\label{d#1}\rm }
\newcommand \bex[1] { \begin{example}\label{e#1}\rm }
\newcommand \bre[1] { \begin{remark}\label{r#1}\rm }

\newcommand \bnota[1] { \begin{notation}\label{n#1}\rm }
\renewcommand {\eth} { \end{theorem} }
\newcommand {\ele} { \end{lemma} }

\newcommand {\epr} { \end{proposition} }
\newcommand {\eco} { \end{corollary} }
\newcommand {\ede} { \end{definition} }
\newcommand {\eex} { \end{example} }
\newcommand {\ere} { \end{remark} }

\newcommand {\enota} { \end{notation} }
\newcommand \thref[1]{Theorem \ref{t#1}}
\newcommand \leref[1]{Lemma \ref{l#1}}
\newcommand \prref[1]{Proposition \ref{p#1}}
\newcommand \coref[1]{Corollary \ref{c#1}}

\newcommand \lb[1]{\label{#1}}

\def \Rset {{\mathbb R}}         
\def \Cset {{\mathbb C}}

\def \tP   {{\widetilde{P}}}

\def \L  {{\mathcal{L}}}           
\def \R  {{\mathcal{R}}}
\def \B  {{\mathcal{B}}}
\def \O  {{\mathcal{O}}}

\def \ga {\gamma}

\def \mt  {\mapsto}
\def \ra  {\rightarrow}           

\def \st  {\ast}                 

\def \ol {\overline}

\def \kk {{\mathrm{k}}}
\def \id { {\mathrm{id}} }
\def \Id { {\mathrm{Id}} }

\def \Ad { {\mathrm{Ad}} }

\def \Lie { {\mathrm{Lie \,}} }
\newcommand{\lag}{\langle}
\newcommand{\rag}{\rangle}
\def \g  {\mathfrak{g}}   
\def \h  {\mathfrak{h}}
\def \k  {\mathfrak{k}}

\def \b  {\mathfrak{b}}
\def \p  {\mathfrak{p}}
\def \r  {\mathfrak{r}}
\def \sl  {\mathfrak{sl}}

\def \rc    {\check{r}}
\def \funG {{\mathcal{C}}^\infty(G)}
\DeclareMathOperator \Span { {\mathrm{Span}} }

\DeclareMathOperator \ad { {\mathrm{ad}} }

\newcommand \bs {\backslash}
\begin{document}
\title
[Triangular Poisson Lie Groups and 
Symplectic Reduction]
{Triangular Poisson Structures on Lie Groups and \\
Symplectic Reduction}
\author{Timothy J. Hodges}
\address{Department of Mathematics, University of Cincinnati, 
Cincinnati, OH 45221-0025}
\email{timothy.hodges@uc.edu}
\thanks{The first author was supported in part by grants from the National 
Security Agency and the Charles
P. Taft Foundation} 
\author{Milen Yakimov}
\address{Department of Mathematics,
University of California,
Santa Barbara, CA 93106}
\email{yakimov@math.ucsb.edu}
\thanks{
The second author was supported in part by NSF grant DMS-0406057
and a junior faculty research incentive grant from UCSB}
\subjclass[2000]{Primary 53D20; Secondary 17B62, 53D17, 17B37}
\date{}

\begin{abstract}
We show that each triangular Poisson Lie group can be decomposed into Poisson 
submanifolds each of which is a quotient of a symplectic manifold. The 
Marsden--Weinstein--Meyer symplectic reduction technique is then used to give
a complete description of the symplectic foliation of all triangular Poisson 
structures on Lie groups. The results are illustrated in detail for
the generalized Jordanian Poisson structures on $SL(n).$ 
\end{abstract}
\maketitle
\section{Introduction}
\lb{1}
A Poisson group structure on a Lie group $G$ is given by 
specifying a Lie bialgebra structure on its Lie algebra $\g$. Such a Lie 
bialgebra structure is said to be triangular if the 
cocommutator is induced from a skew-symmetric solution $r \in \g \wedge \g$
of the classical Yang--Baxter equation. A triangular Poisson Lie group 
is a Lie group equipped with such a Poisson group structure, associated to 
a triangular Lie bialgebra structure on $\g$. If we denote the left and right 
invariant vector fields on $G,$ induced by an element 
$x \in \g \cong T_e G$ by $\L_x$ and $\R_x$ and if 
$r = \sum r^1_i \wedge r^2_i,$
then the Poisson bivector field on $G$ is given by
\begin{equation}
\lb{pi}
\pi_r = \sum \L_{r^1_i} \wedge \L_{r^2_i} - 
\sum R_{r^1_i} \wedge \R_{r^2_i}.
\end{equation}

In this article we give a description of the symplectic leaves of 
all triangular Poisson Lie groups. From the point of view of quantization, 
our results are important because the 
primitive ideals of the quantization are generally closely connected 
with the closures of the symplectic leaves. Since the primitive ideals are 
the annihilators of the irreducible representations of an algebra 
this yields an important link between representations of the 
quantization and symplectic leaves of the Poisson manifold. The most 
well known examples of this phenomenon are the Kostant--Kirillov structure 
on the dual $\g^*$ of a solvable algebraic Lie algebra and the 
standard Poisson structure on a semisimple Lie group. 
In both cases there is a natural bijection between the primitive 
ideals of the quantization and the symplectic leaves of the 
Poisson Lie group. In the former case, the quantization 
is the enveloping algebra and the symplectic leaves are precisely 
the coadjoint orbits \cite{Di}. In the latter case the 
quantization is the Drinfeld--Jimbo quantum group 
which played a major role in many recent developments 
in mathematical physics and pure mathematics \cite{CP}. 
The symplectic leaves in this case 
are closely related to double Bruhat cells
in semisimple Lie groups \cite{Jo, HL1, HKKR, KZ, Y}.
A second important stimulus for studying symplectic leaves 
in Poisson Lie groups is their relation to integrable systems.
Poisson Lie groups provide a very general framework for 
constructing integrable systems \cite{STS3, R, HKKR} 
whose phase spaces are the leaves of those Poisson structures. 
Detailed geometric information for the latter is crucial in 
understanding the properties of such dynamical systems. 

The main tool for studying the leaves of 
a Poisson Lie group $G$ is Semenov-Tian-Shansky's (local) dressing action
of the dual Poisson Lie group $G^*$ on $G.$ 
The problem of explicit 
parametrization of the orbits of this actions is highly nontrivial
and was only solved for specific classes of Poisson Lie groups
see \cite{LS, LW, HL1}, most generally for any Belavin--Drinfeld 
structure on a complex simple Lie group \cite{Y}. The problem 
of understanding the geometry of leaves is open even in the 
simplest case of the standard Poisson structure on a complex
simple Lie group. 
In this paper we show that one can study each triangular Poisson 
structure on a Lie group and in particular its symplectic leaves 
using different methods, namely Marsden--Weinstein--Meyer symplectic 
reduction. Our methods provide an explicit description of the 
symplectic leaves of any triangular Poisson Lie group and provide a
framework for the future study of the geometry of these leaves. 
The approach applies to both complex and real groups. 
All manifolds in this paper can be of either type. 
For shortness we denote
\begin{equation}
\label{k}
\kk = \Cset \quad \mbox{or} \quad \Rset,
\end{equation} 
according to whether the manifold in discussion is complex
or real.

Suppose that $G \times M \to M$ is a free, proper, symplectic action of 
a connected Lie group $G$ on a connected symplectic manifold 
$(M, \pi^{-1} \in \wedge^2 T^* M)$ for a nondegenerate 
Poisson structure $\pi \in \wedge^2 T M.$ 
Any smooth function $f$ on $M$ defines the Hamiltonian
vector field $X_f = df \rfloor \pi$. 
The Lie algebra $\g$ of $G$ acts naturally as vector fields on $M$. Denote 
the vector field on $M$ induced by an element $x \in \g$ by $x_M$.
Recall that a
moment map for this action is a map $J \colon M \to \g^*$ such 
that $X_{J^*(x)} = x_M$. If such a map exists, the action is called {\em weakly 
Hamiltonian}. In general a moment map will not be Poisson or, equivalently, 
$G$ equivariant. It is made such by several twistings 
(equivalently by passing to central extensions): 

1.) The deviation from $G$ equivariance is measured by a one cocycle 
$\B \colon G \to \g^*$ defined by $\B(g) = J(gm) - \Ad^*_g (J(m))$ 
which is independent of
$m \in M$. The coadjoint action of $G$ on $\g^*$ can then be twisted to a new 
geometric action by defining $g*\xi = \Ad^*_g (\xi) + \B(g)$. The moment
map $J \colon M \ra \g^*$ is $G$ equivariant with respect to this new action 
of $G$ on $\g^*.$  

2.) The derivative $\ga = d_e \B,$ which is a Lie algebra 1-cocycle 
$\ga \colon \g \ra \g^*,$ gives rise to the 2-cocycle for $\g$
\begin{equation}
\label{two}
\beta \colon \g \wedge \g \ra \kk, \quad
\beta(x,y) = - \lag \ga(x), y \rag = \lag x, \ga(y) \rag,
\end{equation}
cf. \eqref{k} ($e \in G$ is the identity element).
It can be equivalently defined by
\[
\beta(x, y) = \{ J^*(x), J^*(y) \}
\]
which is necessarily a constant function on $M.$ The cocycle
$\beta$ induces
the affine Poisson structure on $\g^*,$ defined by
\begin{equation}
\label{affine}
\{x,y\} = [x,y] + \beta(x,y), \; \forall x, y \in 
\g \subset \kk[\g^*].
\end{equation}
Alternatively one arrives at it by considering a central extension
of $\g$ by the 2-cocycle \eqref{two} and then restricting the 
Kirillov--Kostant Poisson structure to a hyperplane in its dual, 
identified with $\g^*.$ For this new Poisson structure on $\g^*$ 
the moment map $J \colon M \ra \g^*$ is Poisson. The symplectic
leaves of \eqref{affine} are the orbits of the $\B$-twisted coadjoint 
action of $G$ on $\g^*$. 

Finally the orbit space $M/G$ has a natural manifold structure
and inherits the Poisson structure $p_* (\pi)$ from $M.$ The
Marsden--Weinstein--Meyer \cite{MW, M} reduction theorem states
that the symplectic leaves of $(M/G, p_*(\pi))$ can be described 
as the connected components of $p (J^{-1}(\O))$ where $\O$ is an orbit 
for the $\B$-twisted coadjoint action of $G$ on $\g^*$ and $p$ is the 
natural projection $p \colon M \to M/G$. The symplectic manifolds
$(p (J^{-1}(\O)), p_*(\pi)^{-1})$ are called symplectic 
reductions of $(M, \pi^{-1})$ with respect to $G,$ $\O.$ 
The connectedness
of $G$ and $M$ is essential in establishing the properties 
of the moment map $J$ for the weak Hamiltonian action, 
and the properties of the cocycles $\B,$ $\gamma$ and $\beta.$
On the other hand if these properties are satisfied, then
the symplectic reduction technique is applicable even if
both $G$ and $M$ are disconnected. For more details on weakly 
Hamiltonian actions 
and moment maps we refer to the books \cite{AM, CW}.

Let $(G, \pi_r)$ be a triangular Poisson Lie group. The subspace
$\p = \{(\xi \otimes \Id)(r) \mid \xi \in \g^* \}$ (i.e. the span
of all components of the r-matrix  $r \in \g \wedge \g$) is a Lie
subalgebra of $\g$ which has a canonical structure of a quasi-Frobenius
Lie algebra, see Section \ref{2} for details. Moreover $r$ is a 
non-degenerate skew-symmetric solution of the classical Yang--Baxter equation 
on $\p,$ in particular the map $\rc \colon \p^* \to \p$ is a
linear isomorphism. The inverse map $\rc^{-1} \colon \p \ra \p^*$ is 
a Lie algebra 1-cocycle. Denote by $P$ the connected subgroup 
of $G$ corresponding to $\p.$ (In general $P$ is not a 
closed subgroup of $G.)$ Let $\eta \colon \tP \ra P$ be a connected cover
for which there exists a group cocycle $\B \colon \tP \ra \p^*$
such that $d_e \B = \rc^{-1}.$ In \leref{cosets} we show that all  
double cosets $PtP$ are complete Poisson submanifolds of $G$.
By a complete Poisson submanifold of a given Poisson manifold
we mean a submanifold which is stable under all Hamiltonian flows, 
i.e. a submanifold which is a union of symplectic leaves, 
see Section \ref{3} for details. Thus the symplectic foliation of 
$(G, \pi_r)$ is the union of 
the symplectic foliations of these double cosets.
Each of them can be obtained as a quotient of a  
a symplectic manifold as follows. Fix $t \in G.$ 
The group $\tP \times \tP$ has a natural left invariant 
symplectic structure $(-\L_r, \L_r),$ see \eqref{Lr} and \leref{Invar}
below. In Section \ref{3} we show that the map 
\[
\nu_t \colon \tP \times \tP \ra PtP, \quad
\nu_t(p_1, p_2) = \eta(p_1^{-1}) t \eta(p_2), p_i \in \tP
\]
is Poisson and that the fibers of this map are exactly the left orbits 
of the subgroup 
\[
\tP_t = \{(p_1,p_2) \in \tP \times \tP \mid 
\eta(p_2) = \Ad_t^{-1} \eta_1 (p_1) \}
\]
of $\tP \times \tP.$
Its Lie algebra can be identified with $\p_t = \p \cap \Ad_t(\p).$
We further show that the
action of $\tP_t$ on $\tP \times \tP$ is weakly Hamiltonian
with a moment map
\begin{equation}
\lb{J_t}
J_t(p_1, p_2) = - \B(p_1) + \Ad_{t^{-1}}^* \B(p_2) + \p_t^\perp
\in \g^*/(\p_t^\perp) \cong \p_t^*, 
\end{equation}
where for a subspace $V \subset \g$ by $V^\perp$ we denote its 
orthogonal complement in $\g^*.$ Note that $\B$ takes values in 
$\p^* \cong \g^*/\p^\perp$ and \eqref{J_t} makes sense because
$\p^\perp$ and $\Ad_t(\p^\perp)$ are contained in $\p_t^\perp.$
We refer to \eqref{J_t2} for another way of writing the 
formula for the moment map $J_t(.,.).$ The associated group 
cocycle for this action is $\B_t = J_t|_{\tP_t}$. 
Applying the reduction theory 
described above yields the description of the symplectic leaves of $G$.

\bth{mainthm} Let $\pi_r$ be a triangular Poisson structure on 
a Lie group $G.$
In the above notation the symplectic foliation of $(G, \pi_r)$ 
is the union of the symplectic foliations of all double cosets 
$PtP$. For any $t \in G$ the symplectic leaves of 
$(PtP, \pi_r|_{PtP})$ are the connected components of
$\nu_t(J_t^{-1}(\O))$ where $\O$ is an orbit for the $\B_t$-twisted 
coadjoint action of $\tP_t$ on $\p_t^*$.
\eth

Let us note that the symplectic foliation of any 
Poisson group structure on a complex simple Lie group
can be described on the basis of \thref{mainthm} and 
the results in \cite{Y}, which treated all 
Belavin--Drinfeld Poisson structures (i.e. 
all Poisson structures that are not triangular). 
It would be very interesting 
to find a noncommutative version of our symplectic reduction 
approach to study the primitive spectra of Hopf algebras that 
are quantizations of triangular Poisson structures
on algebraic groups. 

The proof of \thref{mainthm} is given in Section \ref{3}. 
The theorem simplifies in the 
important special case when $r \in \p \wedge \p$ is induced by 
a Frobenius structure 
on $\p.$ This is formulated in \prref{mainthm2}. 
Section \ref{2} contains auxiliary results
on symplectic structures related to triangular Poisson Lie groups
and their symplectic reductions, as well as some 
background on Frobenius and quasi-Frobenius structures on Lie algebras.
Section \ref{4} treats Jordanian Poisson structures 
on $SL(n)$ in detail.   

The results in this paper were obtained in the Fall of 2000 and were 
reported at the conference on Quasiclassical and Quantum Structures
at the Fields Institute (January 2001), and at seminars in Berkeley 
(November 2000) and Cornell (January 2002). We are grateful  
to Alan Weinstein whose comments at one of these talks helped us
to improve our arguments. We are also thankful to him for 
his advice on notation, in particular he suggested to us the term
complete Poisson submanifold. 

The second author would like to thank the organizers of the conference
on Noncommutative Geometry and Representation Theory, Sweden 2004,
for the invitation to participate.
\section{Quasi-Frobenius Lie algebras, symplectic forms and reduction}
\lb{2}
\ble{Invar} 
Let $G$ be a Lie group with Lie algebra $\g$ and
let $r = \sum r^1_i \wedge r^2_i \in \g\wedge \g.$ 
The following left $G$ invariant bivector field
\begin{equation}
\label{Lr}
\L_r = \sum \L_{r^1_i} \wedge \L_{r^2_i}
\end{equation}
on $G$ is Poisson if and only if $r$ is a solution
of the classical Yang--Baxter equation. In this case
\[
\R_r = \sum \R_{r^1_i} \wedge \R_{r^2_i}
\]
is a right invariant Poisson structure on $G.$
\ele

A proof of this Lemma can be found in \cite[Theorem 2.2.2]{CP}.

Fix a Lie group $G$ and a nondegenerate
solution of the classical Yang--Baxter equation 
$r \in \g \wedge \g$ where $\g = \Lie G.$
The nondegeneracy condition means that 
$r \in \g \wedge \g$ defines a nondegenerate
skew-symmetric (i.e. symplectic) form on $\g^*.$
The following relation to quasi-Frobenius
Lie algebras is due to Stolin \cite{S}, see 
for details \cite[Section 3.1.D]{CP} and 
\cite[Section 3.5]{ES}. First the map
\begin{equation}
\label{rcheck} 
\rc \colon \g^* \to \g, \quad \mbox{given by} \quad
\rc(\xi) = (\xi\otimes \Id)r, \; \xi \in \g^*
\end{equation}
is an isomorphism between the underlying Lie algebra $\g^*$ of the dual 
Lie bialgebra of $\g$ and $\g$.

{\em{1. The inverse $\rc^{-1} \colon \g \to \g^*$ is a Lie algebra 
1-cocycle for the coadjoint action of $\g$ on $\g^*$.}}

{\em{2. The bilinear form
\begin{equation}
\lb{beta}
\beta \colon \g \wedge \g \ra \kk, \quad 
\beta(x,y) = - \lag \rc^{-1}(x), y \rag =
\lag x, \rc^{-1}(y) \rag 
\end{equation}
is a nondegenerate 2-cocycle for $\g$ with values in the trivial 
representation, recall \eqref{k}}.}

The second property is the definition of a quasi-Frobenius Lie algebra 
$(\g, \beta)$ and the above construction can be inverted to show that 
the class of those coincides with the class 
of triangular Lie bialgebras $(\g, r)$ with a nondegenerate 
$r$-matrix. \leref{Invar} implies that 
a Lie group $G$ possesses a left invariant symplectic structure 
if and only if $(\Lie(G), \beta)$ is a quasi-Frobenius Lie algebra
for some 2-cocycle $\beta.$ Formulas \eqref{Lr} and \eqref{beta}
give the exact relation between the symplectic form and the 
2-cocycle.

Recall also that a Frobenius Lie algebra is a pair 
$(\g, \xi)$ where $\g$ is a Lie algebra and $\xi \in \g^*$ is 
such that $x \otimes y \ra \lag \xi, [x, y] \rag$ is a 
nondegenerate form on $\g.$ The class of such is the 
subset of all quasi-Frobenius Lie algebras for which
$\beta$ is the coboundary of the 1-cochain 
$x \mt \lag \xi, x \rag,$ or equivalently $\rc^{-1}$ 
is the coboundary of the 0-cochain $\xi \in \g^*.$ 

The following proposition describes under what conditions 
in the above situation the left action of $G$ on $(G, \L_r)$ 
is weakly Hamiltonian. 

\bpr{prop1} Let $G$ be a Lie group for which there exists a
nondegenerate skew-symmetric solution $r \in \g \wedge \g$ 
of the classical Yang--Baxter equation. 
Let $(G,\L_r)$ be the corresponding 
left invariant symplectic structure on $G$. 
Then the following are equivalent.
\begin{enumerate}
\item The left action of $G$ is weakly Hamiltonian.
\item There exists a 1-cocycle $\B \colon G \to \g^*$ such that 
$d\B_e = \rc^{-1}$.
\end{enumerate}
In particular these conditions will hold if $G$ is 
simply connected or if the 1-cocycle 
$\rc^{-1} \colon \g \ra \g^*$ is a coboundary.
\epr

\begin{proof}
If $J \colon G \to \g^*$ is a moment map for the left action of 
$G$ on $G$, then $J - J(e)$ is another moment map which is 
necessarily a 1-cocycle. This follows from the fact that 
for any moment map $J \colon M \ra \g^*$ the function
$\B(g) = J(gm) - \Ad^*_g (J(m))$ does not depend on $m \in M$
and is a 1-cocycle for $G.$ Therefore, to prove the proposition, 
it is enough to observe that a 1-cocycle 
$\B \colon G \to \g^*$ is a moment map if and only if 
$d_e \B = \rc^{-1}$ (where $e$ denotes the identity element 
of $G$).
	
For any $x \in \g$ and $f \in \funG$
\begin{align*}
(X_{\B^*(x)} \cdot f)(e) &= \{ \B^*(x), f(x) \}(e) \\
&= \sum \lag r_i^1,d_e(\B^*(x)) \rag \lag r_i^2, d_e f \rag \\
&= \sum \lag \lag r_i^1,d_e \B(x) \rag r_i^2 , d_e f\rag \\
&= \lag \rc(d_e\B(x)), d_e f \rag.
\end{align*}
On the other hand $x_G = \R_x$ and 
$(\R_x\cdot f)(e) = \langle x, d_e f \rangle$. Hence $X_{\B^*(x),e} = x_{G,e}$ 
if and only if $d_e\B(x)=\rc^{-1}(x)$. This proves the necessity of the 
condition. 

Now assume that $\B$ is a cocycle and that $d_e\B=\rc^{-1}$. 
For $g \in G$ denote 
\[
\ell_g \colon \funG \ra \funG, \quad (\ell_g f)(g_1) = f (g_1 g^{-1}).
\]
Notice that the cocycle condition implies that
\[
\ell_g\B^*(x) = \B^*(\Ad_g x) - \B(g)(x)
\]
and so $X_{\B^*(\Ad_gx)} = X_{\ell_g\B^*(x)}$. Thus
\begin{align*}
(X_{\B^*(x)}\cdot f)(g) &= 
(X_{\ell_{g^{-1}}\B^*(x)} \cdot \ell_{g^{-1}}f)(e) 
= (X_{\B^*(\Ad_{g^{-1}} x)}\cdot \ell_{g^{-1}}f)(e)\\
&= (\R_{\Ad_{g^{-1}}x}\cdot \ell_{g^{-1}} f)(e) =
(\ell_{g^{-1}}\R_x \cdot f)(e)\\
&= (\R_x \cdot f)(g)
\end{align*}
Hence $X_{\B^*(x)} = \R_x=x_G$ and $\B$ is a moment map.
\end{proof}

Suppose that the left action of $G$ on $(G,\L_r)$ is weakly Hamiltonian. 
Let $K$ be a closed subgroup of $G$. The next Proposition 
describes the leaves of the reduction of $(G, \L_r)$
by $K$ in the setting of \prref{prop1}. Let $\k = \Lie K$
and $\pi_{\k^*} \colon \g^\st \ra \k^*$ denote the 
projection associated to the inclusion $\k \subset \g.$
If $\B$ is a moment map for the action of $G$ on itself which is a
group cocycle, then $J_K =\pi_{\k^*} \circ \B \colon G \to \k^*$ is a 
moment map for the left action of $K$ on $G$. The cocycle for $K$ 
associated to this weakly Hamiltonian map is $\B_K = \pi_{\k^*}\circ \B|_K$
(because $\B(e)=0)$. Let $\nu \colon G \to K \bs G$ be the natural projection.
It induces the Poisson structure $\nu_*(\L_r)$ on 
$K\bs G$ for which the symplectic reduction technique 
yields the following description 
of its symplectic leaves.

\bpr{cor1} Let $G$ be a Lie group for which 
$(\g = \Lie G, \beta)$ is a Frobenius Lie algebra
and let $r \in \g \wedge \g$ be the corresponding nondegenerate
$r$-matrix. Assume that there exists a group 1-cocycle 
$\B \colon G \ra \g^*$ with 
$d_e \B = \rc^{-1}.$ If $K$ is a closed subgroup of $G$ then  
the symplectic leaves of $(K \bs G, \nu_*(\L_r))$ 
are the connected components of $\nu(J_K^{-1}(\O))$ where $\O$ 
is an orbit for the $\B_K$-twisted coadjoint action of $K$ on $\k^*$.
\epr

Recall that the $\B_K$-twisted coadjoint action of $K$ on $\k^*$
is the action $k * \xi = \Ad^*_k(\xi) + \B_K(k).$ 

\prref{cor1} simplifies in the case of Frobenius Lie algebras.

\bpr{cor2} Let $G$ be a Lie group for which 
$(\g = \Lie G, \xi \in \g^*)$ 
is a quasi-Frobenius Lie algebra and let 
$r \in \g \wedge \g$ be the corresponding nondegenerate
$r$-matrix. Then the symplectic leaves of $(K \bs G, \nu_*(\L_r))$ 
are the connected components of $\nu(I_K^{-1}(\O))$
where 
\[
I_K \colon G \ra \k^* \; \; \mbox{is given by} \; \;
I_K(g)= \pi_{\k^*}\Ad^*_g(\xi), g \in G
\]  
and $\O$ is a coadjoint orbit of $K.$ 

Equivalently, the symplectic leaf of $(K \bs G, \nu_*(\L_r))$
through $Kg$ is the connected component of the submanifold
\[
X_g = \{ Kg' \in K \bs G \mid g'\in G, \;
\Ad^*_{g'} (\xi) \in \Ad^*_g(\xi) + \k^\perp\} \subset K \bs G,
\]  
containing $Kg.$
\epr

It is clear that only the coadjoint orbits in 
$\pi_{\k^*}(\Ad_G^*(\xi))$ 
have nontrivial preimages and thus lead to symplectic leaves.

\begin{proof} In this case $\rc^{-1} \colon \g \ra \g^*$ is coboundary:
$\rc^{-1}(x) = \ad_x^*(\xi)$ for all $x \in \g.$ It integrates to 
the group coboundary $\B \colon G \ra \g^*,$ 
$\B(g) = \Ad_g^*(\xi) - \xi.$ 
Therefore the moment map $J_K$ is 
\[
J_K(g) = \pi_{\k^*}\B(g) = \pi_{\k^*} \Ad_g^*(\xi) - 
\pi_{\k^*}(\xi) = I_K(g) - \pi_{\k^*}(\xi).
\]
The corresponding group cocycle for $K$ is the coboundary
of $\pi_{\k^*}(\xi)$ and the related twisted coadjoint orbits
are just the translation of the coadjoint orbits of $K$ by
$-\pi_{\k^*}(\xi).$ This implies the first part of the proposition.

For the second part one first observes that the symplectic leaf of 
$(K \bs G, \nu_*(\L_r))$ through $Kg$ is the connected component of 
\[
\{ Kg' \in K \bs G \mid g'\in G, \;
\Ad^*_{g'} (\xi) \in \Ad^*_K \Ad^*_g(\xi) + \k^\perp\},
\]  
containing $Kg.$ This is exactly the set $X_g$ because $\k^\perp$ is 
$\Ad_K^*$ stable.
\end{proof}
\section{Symplectic leaves of triangular Poisson Lie groups}
\lb{3}
Throughout this section we fix a triangular Poisson structure
$\pi_r$ on a Lie group $G,$ cf. \eqref{pi}, where
$r \in \g \wedge \g$ is a skew-symmetric 
solution of the classical Yang--Baxter equation
on $\g = \Lie G.$ The subspace
\begin{equation}
\lb{p}
\p = \Span \{(\xi \otimes \Id) r \mid \xi \in \g^* \}
\end{equation}
is a Lie subalgebra of $\g$ because it is the image of the 
Lie algebra homomorphism $\rc \colon \g^* \ra \g,$
$\rc(\xi) = (\xi \otimes \Id) r$ from the underlying Lie algebra 
of the dual Lie bialgebra $\g^*$ of $\g$ to $\g,$ see e.g. 
\cite{ES}.
It is clear that $r \in \p \wedge \p$ is a non-degenerate
solution of the classical Yang--Baxter equation on $\p$
which turns $(\p, \beta)$ into a quasi-Frobenius Lie algebra
for $\beta \colon \p \wedge \p \ra \kk$ given by \eqref{beta},
recall \eqref{k}. 
Denote the connected subgroup of $G$ with tangent Lie algebra 
$\p$ by $P.$ Note that in general $P$ is not a closed subgroup 
of $G.$ Recall from Section 2 that $\rc \colon \p^* \ra \p$
is a linear isomorphism and its inverse 
$\rc^{-1} \colon \p \ra \p^*$ is a 1-cocycle. It 
integrates to a group cocycle
\begin{equation}
\lb{B}
\B \colon \tP \ra \p^*, \; d_e \B = \rc^{-1}
\end{equation}
for an appropriate (connected) cover $\tP$ of $P.$ 
Denote the projection $\eta \colon \tP \ra P.$
We can always choose $\tP$ to be the simply connected 
group with Lie algebra $\p.$ If $\rc^{-1}$ is a 
coboundary of $\xi \in \p^*$ we can take 
$\tP = P$ and set $\B(p) = \Ad^*_p(\xi) -\xi,$
$\forall p \in P.$    

First we define the notion of a complete Poisson submanifold. 
\bde{ps} 
Let $(M, \pi)$ be a Poisson manifold. A (not necessarily closed) 
submanifold $X$ of $M$ is called a complete Poisson submanifold if $X$ is 
stable under all Hamiltonian flows on $M.$
\ede

If $X$ is a complete Poisson submanifold of $(M, \pi)$ then 
\begin{equation}
\lb{rest}
\pi_x \in \wedge^2 T_x X \; \forall x \in X,
\end{equation}
i.e. $X$ is a Poisson submanifold of $(M, \pi)$ in the 
local sense, see e.g. \cite{V}.
Conversely a Poisson submanifold $X$ is not necessarily a complete 
Poisson submanifold of $(M, \pi),$ this is only true if $X$ is a closed
submanifold of $M.$ A submanifold $X$ of $(M, \pi)$ is a complete Poisson 
submanifold if and only if it is a union of symplectic
leaves of $M.$ The following Lemma clarifies 
the passage from Poisson submanifolds to complete
Poisson submanifolds.

\ble{pass} Assume that $(M, \pi)$ is a Poisson manifold which can be 
decomposed into the disjoint union of (not necessarily closed) 
submanifolds $M = \sqcup_{\alpha \in A} X_\alpha$ satisfying
\begin{equation}
\lb{rest1}
\pi_x \in \wedge^2 T_x X_\alpha \; \forall x \in X_\alpha,
\end{equation}
i.e. all $X_\alpha$ are Poisson submanifolds of $(M, \pi).$
Then all $X_\alpha$ are complete Poisson submanifolds of 
$(M, \pi).$
\ele

\begin{proof} If $\gamma \colon (a, b) \ra M$ is a Hamiltonian path 
then \eqref{rest1} implies that $\gamma^{-1}(X_\alpha)$ is an open subset 
of $(a,b)$ for all $\alpha \in A.$ Since $(a,b)$ is connected the image of 
$\gamma$ sits within a single submanifold $X_\alpha$ of $(M, \pi)$
which is easily seen to complete the proof of the Lemma.
\end{proof} 

\ble{cosets} 
Assume that $(G, \pi_r)$ is a triangular Poisson Lie 
group and $P$ is the connected subgroup of $G$ whose Lie algebra is 
the span of all components of $r \in \Lie G \wedge \Lie G,$ recall 
\eqref{p}. Then all double cosets $P t P,$ $t \in G$ are complete 
Poisson submanifolds of $(G, \pi_r).$ 
\ele

\begin{proof}
We apply \leref{pass} to the decomposition
\[
G = \sqcup_{[t] \in P \bs G/P} P t P
\]
over some representatives of all $(P,P)$-double cosets of $G.$ 
Observe that 
\[
\pi_r(x) \in \wedge^2 T_x (PtP) \; \forall t \in G, x \in PtP
\]
because of the definition \eqref{pi} of $\pi_r$ 
and the fact that the left and right-invariant
bivector fields $\L_r$ and $\R_r$ on $G$ are tangent to $PtP$ for 
all $t \in G.$  
\end{proof}
 
In the rest of this section we describe how each of the Poisson submanifolds
$PtP$ of $(G, \pi_r)$ is obtained as a quotient of a symplectic manifold in the 
framework of \prref{cor1}. In particular
we complete the proof of \thref{mainthm}. Consider the left-invariant 
Poisson structure $(-\L_r, \L_r)$ on $\tP\times \tP$ coming from the 
skew-symmetric solution $(-r,r)$ of the classical Yang--Baxter on 
$\p \oplus \p.$ Recall from the Introduction the definition of the subgroup
\begin{equation}
\lb{p_t}
\tP_t = \{ (p_1,p_2) \in \tP \times \tP \mid 
\eta(p_2) = \Ad_t \eta(p_1)\} \subset \tP \times \tP.
\end{equation}
It is clear that $\tP_t$ is a closed subgroup 
of $\tP \times \tP.$ Observe that the fibers of the projection 
\begin{equation}
\lb{nu}
\nu_t \colon \tP \times \tP \ra P t P, \quad 
\nu_t(p_1, p_2) = \eta(p_1^{-1}) t \eta(p_2), \; p_i \in \tP  
\end{equation}
are exactly the left $\tP_t$ orbits on $\tP \times \tP.$ 
It is easy to see that 
\[
\nu_{t *} (- \L_r, \L_r) = - \R_r + \L_r = \pi_r
\]
using the fact that for any triangular Poisson Lie group $(G, \pi_r)$
we have an isomorphism of Poisson manifolds 
$(G, \L_r) \ra (G, \R_r)$ given by $g \mt g^{-1}.$ Thus
$(PtP, \pi_r|_{PtP})$ is the quotient of 
$(\tP \times \tP, (-\L_r, \L_r))$ by the left action of 
the subgroup $\tP_t$ of $\tP \times \tP.$ To apply \prref{cor2}
we need the explicit form of the corresponding moment maps.
Let
\[
\p_t = \p \cap \Ad_t(\p).
\]
We will identify the Lie algebra of $\tP$ with $\p_t$ by the 
isomorphism
\[
\Lie(\tP_t) = \{ (x, \Ad_t(x) ) \mid x \in \p_t \} \cong \p_t
\]
given by the projection into the first factor.

Recall that $\B \colon \tP \ra \p^*$ is a group 1-cocycle 
with $d_e \B = \rc^{-1}.$ Thus the map 
$(p_1, p_2) \mt (-\B(p_1), \B(p_2)) \in \p^* \oplus \p^*$
is a 1-cocycle for $\tP \times \tP$ and  
according to \prref{prop1} it is a moment map for the left 
action of $\tP \times \tP$ on $(\tP \times \tP, (-\L_r, \L_r)).$ 
If 
\begin{equation}
\label{eq}
\pi'_t \colon \p^* \ra \p_t^* = ( \p \cap \Ad_t(\p))^* 
\; \mbox{and} \; \pi''_t \colon \p^* \ra (\Ad_t^{-1} \p_t)^* 
= ( \Ad_t^{-1}(\p) \cap \p)^* 
\end{equation}
denote the natural projections, then the map 
$(p_1, p_2) \mt (-\pi'_t \B(p_1), \pi''_t \B(p_2)) \in 
\p_t^* \oplus (\Ad_t^{-1} \p_t)^*$
is a moment map for the left action of 
$\eta^{-1}(P \cap \Ad_t P) \times \eta^{-1} (\Ad_t^{-1} P \cap P)$ 
on $(\tP \times \tP, (-\L_r, \L_r)).$ This finally means that 
the map 
\begin{equation}
\lb{J_t2}
J_t(p_1, p_2) = - \pi'_t \B(p_1)+ \Ad_{t^{-1}}^* \pi''_t \B(p_2), 
\quad 
J_t \colon \tP \times \tP \ra \p_t^*
\end{equation}
is a moment map for the left action of $\tP_t$
on $(\tP \times \tP, (-\L_r, \L_r)).$ Eq. \eqref{J_t2} 
is just another way of 
writing the moment map $J_t(.,.)$ from \eqref{J_t}. 
Now \thref{mainthm} follows from \prref{cor1}.

\thref{mainthm} simplifies in two special cases.

\bco{qFr} Assume that $\pi_r$ is a triangular Poisson structure 
on a Lie group $P$ such that $r \in \p \wedge \p$ is nondegenerate
(i.e. $r$ corresponds to a quasi-Frobenius structure on $\p$). If
the 1-cocycle $\rc^{-1} \colon \p \ra \p^*$ lifts to a group 
cocycle $\B \colon P \ra \p^*$, then the symplectic leaves
of $(P, \pi_r)$ are the connected components of $\B^{-1}(\O)$
for $\O$ a coadjoint orbit in $\p^*.$
\eco 

To deduce \coref{qFr} from \thref{mainthm} observe that 
$P = P e P,$ $P_e$ is the diagonal of $P \times P,$
and $\B_t \colon P_e \ra \p^*$ vanishes
because $\B_t(p,p) = - \B(p) + \B(p)=0$ for all $p \in P.$
In addition
\[
J_e(p_1, p_2) = -\B(p_1) + \B(p_2) = 
\Ad_{p_1}^*(\B(p_1^{-1}) + \Ad_{p_1^{-1}}^* \B(p_2) )
=\Ad_{p_1}^* \B(p_1^{-1}p_2),
\]
for all pairs $(p_1, p_2) \in P \times P.$

\coref{qFr} was independently obtained by Diatta and Medina in \cite{DM}.

Secondly \thref{mainthm} also simplifies in the case when 
$\r \in \g \wedge \g$ is induced by a Frobenius structure on 
a subalgebra of $\g.$ We will omit the proof 
of the following result which 
is easily obtained from \prref{cor2}.

\bpr{mainthm2} Assume that $(G, \pi_r)$ is a triangular Poisson Lie group
for which $r \in \p \wedge \p,$ recall \eqref{p},
is induces by a Frobenius structure $\xi \in \p^*.$
Then the symplectic leaves of each Poisson submanifold
$PtP,$ $(\forall t \in G)$ are the connected components 
of $\nu_t(I_t^{-1}(\O))$ where $\O$ is a coadjoint orbit of 
$\tP_t$ and  
\[
I_t \colon \tP \times\tP \ra \p_t^* \cong \g^*/\p_t^\perp
\; \; \mbox{is given by} \; \;
I_t(p_1,p_2)= - \Ad^*_{p_1}(\ol{\xi})
+ \Ad^*_{t^{-1} p_2}(\ol{\xi}) + \p_t^\perp, p_i \in \tP.
\]  
Here $\ol{\xi}$ is an arbitrary preimage of $\xi$ under the projection
$\g^* \ra \p^*.$ Only the coadjoint orbits in 
$\pi'_t(\Ad_P^*(\xi)) +  \Ad_{t^{-1}}^* \pi''_t(\Ad_P^*(\xi))$ 
have nontrivial preimages, cf. \eqref{eq}.

Equivalently, the symplectic leaf of $(PgP, \pi_r|_{PgP})$
through $PgP$ is the connected component of the submanifold
\[
Y_g = \{ p_1^{-1} g p_2 \mid p_i\in \tP, \;
(\Ad^*_{p_1}(\ol{\xi}) -  \ol{\xi}) -
\Ad^*_{t^{-1}}(\Ad^*_{p_2}(\ol{\xi}) - \ol{\xi})
\in \p_t^\perp
\},
\]
containing $PgP.$
\epr
\section{Jordanian Poisson structures}
\lb{4}
Some of the most interesting examples of triangular Poisson structures
come from Frobenius parabolic subalgebras of complex simple Lie algebras,
see \cite{GG}.

Assume that $G$ is a complex simple Lie group and that $P$ is a 
parabolic subgroup whose Lie algebra $\p$ has a Frobenius structure 
$\xi_0 \in \p^*.$ It gives rise to a nondegenerate $r$-matrix 
$r \in \p \wedge \p \subset \g \wedge \g,$ see Section \ref{2}.
In this section we study the corresponding Poisson structure 
$\pi_r$ on $G.$

Denote by $W_L$ the Weyl group of the pair $(L,T)$
where $L$ is a Levi factor of $P$ and $T$ is a maximal torus of $L.$ 
Denote by $W$ the Weyl group of the pair $(G,T).$
The Bruhat decomposition of $G$ with respect to the 
parabolic subgroup $P$ is
\[
G = \bigsqcup_{w \in W_L\backslash W / W_L} PwP.
\]
Recall that every double coset $W_L\backslash W / W_L$
has a unique element of a minimal length, see e.g. 
\cite{Car}. 
The set of such is usually denoted by ${}^{W_L} W^{W_L}$
and is identified with $W_L\backslash W / W_L.$ \prref{mainthm2}
now implies the following description of the symplectic leaves
of $(G, \pi_r).$

\bco{par} In the above setting the Poisson Lie group 
$(G, \pi_r)$ decomposes into the complete 
Poisson submanifolds $PwP$
for $w \in {}^{W_L} W^{W_L}.$
If $\nu_w \colon P \times P \ra P w P$ denotes the 
projection $\nu_w(p_1, p_2) = p_1^{-1} w p_2$ then   
the symplectic leaves of $PwP,$ are the connected components
of $\nu_t(I_t^{-1}(\O))$ where $\O$ is a coadjoint orbit of
$\tP_t$ and
\[
I_t \colon \tP \times\tP \ra \p_t^* \cong \g^*/\p_t^\perp
\; \; \mbox{is given by} \; \;
I_t(p_1,p_2)= - \Ad^*_{p_1}(\ol{\xi})
+ \Ad^*_{t^{-1} p_2}(\ol{\xi}) + \p_t^\perp, p_i \in \tP.
\]
Here $\ol{\xi}$ is an arbitrary preimage of $\xi$ under the projection
$\g^* \ra \p^*.$
\eco

The simplest examples of this type are Gerstenhaber--Giaquinto's
\cite{GG} generalized Jordanian Poisson structures on $SL(n),$ 
induced from the $r$-matrices 
\[
r = n \sum_{i<j} \sum_{k=i}^{j-1} E_{k,i}\wedge E_{i+j-k-1,j} + 
\sum_{i,j}(j-1) E_{j-1,j} \wedge E_{i,i}.
\]
Such an $r$-matrix corresponds to the maximal parabolic 
subalgebra of $\sl(n),$ associated to an extremal root
\[
\p = \Span \{ E_{i,j} \mid i <n, j \leq n \},
\]
and is induced by the following Frobenius structure on $\p,$
\cite[Theorem 5.9]{GG} 
\[
\xi = -n \sum_{i=1}^{n-1} E_{i \; i+1}^*.
\]
Here $\{E_{ij}^*\}_{i,j=1}^n$ denotes the dual basis to the basis
$\{E_{ij}\}_{i,j=1}^n$ of the space of $n \times n$
matrices. 
The Weyl group of $SL(n)$ is identified with $S_n$ and 
its subgroup corresponding to the Levi subalgebra of $\p,$ 
containing the Cartan subalgebra of diagonal matrices, is 
$S_{n-1}.$ The related minimal length elements, parametrizing 
the Poisson submanifolds in \coref{par}, are
\[
{}^{S_{n-1}} S_n^{S_{n-1}}= \{\id, (n-1 \; n) \}.
\]
The subalgebra $\p_{(n-1 \; n)}$ is easily described:
\[
\p_{(n-1 \; n)} = \Span \{ E_{i,j} \mid i <n, j \leq n \}.
\]
It does not depend on the representative of $(n-1 \; n) \in S_n$ 
in the normalizer of $T.$

In the rest of this section we provide full details for the case $n=2.$ 
We will use the standard basis $\{E,F, H\}$ of $\sl(2)$
and the corresponding dual basis $\{E^*, F^*, H^*\}$ of 
$\sl(2)^*.$ The Jordanian Poisson structure 
on $G=SL(2)$ corresponds to the $r$-matrix $r = E \wedge H$. 
The latter is induced from the Borel subalgebra $\b$ 
generated by $E$ and $H.$ If we identify $E^*$ and $H^*$ with their 
images in $\b^*,$ then the Frobenius structure on $\b^*$ related to $r$ 
is simply
\[
\xi = - E^*.
\]
The cocycle $\rc^{-1} \colon \b \ra \b^*$ integrates to the group cocycle
\[
\B(g)= E^* - \Ad_g^*(E^*)= (1 - a^{-2}) E^* -2a^{-1}b H^* \in \b^*, 
\quad
g=\begin{bmatrix} a & b\\0&a^{-1}\end{bmatrix} \in B.
\]
where $B$ denotes the Borel subgroup of $SL(2),$ 
consisting of upper triangular matrices. The decompositions 
of $(SL(2), \pi_r)$ into complete Poisson submanifolds from 
\coref{par} is just the Bruhat decomposition 
$G = B \cup BsB.$ Here 
$s$ is a representative of the nontrivial 
element of the Weyl group in the normalizer of the 
maximal torus $T$ of diagonal matrices, say 
$s = \begin{bmatrix} 0 & 1\\-1&0\end{bmatrix}$.

The symplectic leaves of $(SL(2), \pi_r)$ 
are then of two types -- ``$s$'' and ``$e$''. The 
leaves of type $e$ are the inverse images $\B^{-1}(\O)$ 
for $\O$ a coadjoint orbit of $\b^*$. The coadjoint orbits of $\b^*$ 
are the points $\alpha H^*$ and the complement to the set of these points. 
The former give rise to the symplectic points
\[
\left\{ \pm\begin{bmatrix} 1 & - \alpha/2\\0&1\end{bmatrix} \right\}.
\]
If $S$ denotes the set of such points, the inverse image 
of the 2 dimensional coadjoint orbit in $\b^*$ is the complement 
of $S$ in $B$ which is a single symplectic leaf.

Now consider the case when $t=s$. 
Then $\b_s=\b \cap \Ad_s(\b) = \h= \Lie T$
and $I_s \colon B \times B \to \h$ is given by
\[
I_s(g_1,g_2) = \Ad_{g_1}^*(E^*) - \Ad_{s g_2}^*(E^*) + \h^\perp = 
2(a_1^{-1}b_1+a_2^{-1}b_2)H^*, \quad
g_i=\begin{bmatrix} a_i & b_i\\0&a_i^{-1}\end{bmatrix}.
\]
The coadjoint orbits of $G_s=T$ are just the points $\alpha H^*.$
Thus the symplectic leaves of $BsB$ are the connected components
of
\[
\nu_s(I_s^{-1}(\alpha H^*)) =
\{ g_1^{-1} s g_2 \mid 2(a_1^{-1}b_1 + a_2^{-1}b_2)=\alpha \}.
\]
The latter is easily seen to be the locally closed variety 
$\mathcal{V}(2(T_{11}-T_{22}) +\alpha T_{21})- B$ which is
connected. By $T_{ij}$ we denote the standard coordinate 
functions on $SL(n).$

Thus from an algebro-geometric point of view, the leaves can be described 
as follows. The zero dimensional symplectic leaves are the points of the 
variety $Z = \mathcal{V}(T_{11}^2-1, T_{21})$. The two dimensional leaves 
are the locally closed varieties 
$L_{[\lambda,\mu]} = \mathcal{V}(\lambda(T_{11}-T_{22}) -\mu T_{21})- Z$ 
where $ [\lambda, \mu] \in \mathbb{P}^1(\mathbb{C})$. Note that the 
closures of the two dimensional leaves all contain $Z$ and that the 
intersection of the closures of any two distinct two-dimensional leaves is 
precisely $Z$.

\end{document}